\documentclass[11pt]{article}
\title{Short Proof of Dirichlet's Principle}
\author{H. N. Friedel}
\usepackage{amsmath}
\usepackage{latexsym}
\usepackage{graphicx}    
\begin{document}
\maketitle
\begin{abstract}
A standard Hilbert-space proof of Dirichlet's principle is simplified, using an observation that a certain form of min-problem has unique solution, at a specified point.
This solves Dirichlet's problem, after it is recast in the required form (using the Poincar\'e/Friedrichs bound and Riesz representation).
The solution's dependence on data is linear and continuous; and the solution is invariant under certain changes of data, away from the border of the region where Dirichlet's problem is given.
If that region is regular enough for functions on it to have border-traces, then the problem can be stated and solved in terms of border-data.
\end{abstract}
\section{Introduction}
Recall that Riemann's impressive use (in 1851, and later) of the unproved Dirichlet principle provoked a search for proof.
The search grew intense after Weierstrass showed (1870) that a min-problem, such as Dirichlet's, may lack solution.
The difficulty of the challenge is suggested by the eminence of the participants, and by the long wait until proofs finally came: Poincar\'e (1887), Hilbert (1899).
The subject's development (with references) is surveyed in [G].

Natural follow-up tasks include generalization of the principle (such as by weakening premises) and simplification of its proofs.
Several forms of the principle arose (as in [G], [Z] \textit{et al}); any two of these are closely related, though often not strictly equivalent.
This note simplifies a modern proof (in [Z]) of a standard general form of the principle.

To prepare to state Dirichlet's problem in the form treated here, let $\Omega$ denote a non-empty bound open part of $\mathbf{R}^n \,$.
Early treatments assumed some regularity for its border, $\mathrm{bdr}(\Omega)$, often making it a hypersurface in $\mathbf{R}^n \,$, suitable for assignment of data; we defer such assumptions.
As usual, $\mathcal{L}_2 (\Omega)$, $W_2^1 (\Omega)$, $\mathring{W}_2^1 (\Omega)$ denote (respectively) the Lebesgue space of square-summable functions (modulo null measure, on $\Omega$), the Sobolev space of functions with square-summable weak-rates, and its subspace of functions vanishing on the border.

Each of these three spaces has a scalar-product, or ``bracket'', that makes it a Hilbert space.  Denote the bracket on $\mathcal{L}_2$ by $(\cdot|\cdot)_2 \,$.
To help denote brackets for the other two spaces, first define ``grad-bracket'' on $W_2^1 \,$:
$$ (u|v)_{\nabla} \; = \; \int_{\Omega}(\nabla u \, | \, \nabla v) \; = \; \int_{\Omega} \sum \partial_i u \, \partial_i v \ .$$
Naturally, grad-bracket is so-called because it uses gradient, $\nabla u \, = \, (\partial_1 u, \, \dots \, \partial_n u)$.
It yields a seminorm: $(u|u)_{\nabla} \, = \, \|u\|_{\nabla}^2 \,$.
On $W_2^1 \,$, grad-bracket is not a true bracket ($\|c\|_{\nabla} \, = \, 0$ for constant $c$), but we get a true bracket by adding the square-sum bracket:
$ (\cdot|\cdot)_{1,2} \; = \; (\cdot|\cdot)_2 \; + \; (\cdot|\cdot)_{\nabla} \; .$
On the subspace $\mathring{W}^1_2 \,$, grad-bracket suffices (without adding $(\cdot|\cdot)_2 \,$),
thanks to the Poincar\'e/Friedrichs bound:
\begin{equation}
\|u\|_2 \; \, \leq \; \, a \, \|u\|_{\nabla} \ \ \; \mathrm{if} \ \; u \in \mathring{W}_2^1 \; . \label{Poin/F}
\end{equation}
Here ($a>0$) depends only on region $\Omega$.  (\ref{Poin/F}) ensures grad-norm $\|\cdot\|_{\nabla}$ is a norm on $\mathring{W}_2^1 \,$, equivalent to original norm $\|\cdot\|_{1,2} \;$.

It will pay to re-formulate the classical notion of border-data, $\{ g_0 \,:\mathrm{bdr}(\Omega) \to \mathbf{R} \}$, as follows.
First, view $g_0$ as ``trace'', or restriction-to-border, of natural extension $\bar{u}$ of sought function $u$ to closure $\bar{\Omega}$.
(A simple case has $\bar{u}$ continuous.  Recall trace exists if $\Omega$ and $u$ behave well-enough.)
Express this view as `` $u$ restricts to $g_0$ '', or ``$g_0$ restricts $u$''.
Next, view $\{ g_0 \,, u \}$ the other way: ``$g_0$ extends to $u$'', or ``$u$ extends $g_0 \,$''; this means $g_0$ extends to $\bar{u}$  (on $\bar{\Omega}$) and $\bar{u}$ restricts to $u$ (on $\Omega$).
Call $u$ a ``minimal'' extension of $g_0 \,$, if $u$ solves Dirichlet's min-problem.
Certainly, if $g_0$ is to have a minimal extension, it must have some extension; a mild premise is that $g_0$ does extend to some $g \in W_2^1 (\Omega)$.
Dirichlet's principle promises a unique minimal extension for border-data with some extension.

After we fix an ``initial extension'' $g$ of $g_0 \,$, then $g$ includes information about $g_0 \,$, and statements with $g_0$ may be recast in-terms-of $g$:
in particular, for the sought function $u \in W^1_2 $, express ``$u$ extends $g_0$'' as ``$u - g \in \mathring{W}_2^1 \, $'', i.e. on the border we have ($u-g=0$).
Thus we ``hide'' $g_0 \,$, and present $g$ instead.
This way is more general; it works even when region $\Omega$ is not regular enough to ensure existence of border-traces:
view each equi-class ($g+ \mathring{W}^1_2$) in quotient-space $ \, W^1_2 \Big/ \mathring{W}^1_2 \,$ as a ``virtual border-function'';
and view $g$ (and every member of its equi-class) as ``extending'' that virtual border-function.
Thus we expect the sought minimal-extension $u$ of $g_0$ should not depend on particular choice $g$ for initial-extension of $g_0 \,$; $u$ should depend on equi-class ($g+ \mathring{W}^1_2$).
Dirichlet's principle says each equi-class in $ \, W^1_2 \Big/ \mathring{W}^1_2 \,$ has a unique minimal member.

Now we are prepared to state Dirichlet's min-problem, and its ``critical equation'' (that seeks a critical-point of the function to minimize). \\
\textbf{Dirichlet's Problem.}  Seek $u \in W_2^1 \,$, given $f \in \mathcal{L}_2 $ and $g \in W_2^1 \, $, with
\begin{equation}
\frac{1}{2} \, \|u\|_{\nabla}^2 \; - \; (f|u)_2 \; = \; \min! \; , \quad u-g \in \mathring{W}_2^1 \; . \qquad \Box  \label{Diri-prob}
\end{equation}
The following critical-equation for (\ref{Diri-prob}) also seeks $u \in W_2^1 \,$, given $f \in \mathcal{L}_2 $ and $g \in W_2^1 \,$.
\begin{equation}
(u|\phi)_{\nabla} \; = \; (f|\phi)_2 \ \; \mathrm{if} \ \phi \in \mathring{W}_2^1 \; \ \; \mathrm{and} \ u-g \in \mathring{W}_2^1 \; . \label{Diri-eqn}
\end{equation}
\textbf{Dirichlet's Principle.}  Each of \{(\ref{Diri-prob}), (\ref{Diri-eqn})\} has unique solution, and the two solutions are equal.
\section{Short Proof}
Recall the following simple fact.  It is the key to a short proof of the principle. \\
\textbf{Lemma.} For any real Hilbert space $X$, a min-problem (seek $x \in X$, given $p \in X$) of the form
\begin{equation}
\frac{1}{2} \, \|x\|^2 \; \, - \; ( p \, | \, x ) \ \; = \  \min ! \label{bra-min}
\end{equation}
has unique solution ($x=p$). \\
\textbf{Proof.}  Adding a constant to a function doesn't change its min-points; in (\ref{bra-min}), add $\; \frac{1}{2} \|p\|^2 \;$:
$$ \frac{1}{2} \, \|x\|^2 \; - \; ( p \, | \, x ) \; + \;  \frac{1}{2} \, \|p\|^2 \; = \; \frac{1}{2} \, ( x-p \, | \, x-p ) \; = \; \frac{1}{2} \, \|x-p\|^2 \; .$$
The preceding function (of $x$) has unique min-point: $x=p$.  \textbf{Done.} \\
(\ref{bra-min}) is stated in terms of bracket, so call (\ref{bra-min}) a ``bracket-min problem''; the ``bracket-min lemma'' solves it.

\textbf{Proof} of Dirichlet's Principle.  Put $u-g=v \in \mathring{W}_2^1 \,$.
Adding a constant to a function preserves min-points; so (\ref{Diri-prob}) is equivalent to the following min-problem. \\
Seek $v \in \mathring{W}_2^1 \,$, given $f \in \mathcal{L}_2 $ and $g \in W_2^1 \,$, with
\begin{equation}
\frac{1}{2} \, \|v\|_{\nabla}^2 \; - \; (f|v)_2 \; + \; (g|v)_{\nabla} \; = \; \min! \; . \label{Diri-prob2}
\end{equation}
Let $\Lambda$ denote the linear function in (\ref{Diri-prob2}):  $\{ \, \mathring{W}_2^1 \owns \phi \mapsto (f|\phi)_2 \, - \, (g|\phi)_{\nabla} \, \}$.
Use (\ref{Poin/F}) to bound $\Lambda$:
\begin{equation}
|\Lambda \phi| \; \leq \; | \, (f|\phi)_2 \, | \; + \; |(g|\phi)_{\nabla}| \; \leq \; \|f\|_2 \, \|\phi\|_2 \; + \; \|g\|_{\nabla} \, \|\phi\|_{\nabla} \; \leq
\; \big( a \, \|f\|_2 \, + \, \|g\|_{\nabla} \big) \, \|\phi\|_{\nabla} \quad \mathrm{if} \ \; \phi  \in \mathring{W}_2^1 \; . \label{Diri-prin}
\end{equation}
Since $\Lambda$ has bound, Riesz-representation gives a unique point $p \in \mathring{W}_2^1$ that represents $\Lambda$:
\begin{equation}
(p|\phi)_{\nabla} \; = \; \Lambda \phi \; = \; (f|\phi)_2 \, - \, (g|\phi)_{\nabla} \quad \mathrm{if} \ \phi  \in \mathring{W}_2^1 \; . \label{Diri-prob4}
\end{equation}
Hence (\ref{Diri-prob2}) can be written as a bracket-min problem: seek $v \in \mathring{W}_2^1$, given $p \in \mathring{W}_2^1 \,$, with
\begin{equation}
\frac{1}{2} \, \|v\|_{\nabla}^2 \; - \; (p|v)_{\nabla} \; = \; \min! \; . \label{Diri-prob3}
\end{equation}
(\ref{Diri-prob3}) has unique solution ($v=p$), by the bracket-min lemma; then (\ref{Diri-prob2}) has the same; and (\ref{Diri-prob}) has unique solution ($u=p+g$).
(\ref{Diri-prob4}) implies this point ($u=p+g$) also solves (\ref{Diri-eqn}).
Solution of (\ref{Diri-eqn}) is unique: if $\tilde{u}$ solves (\ref{Diri-eqn}), then $ u - \tilde{u} = (u - g) - (\tilde{u} - g) \in \mathring{W}_2^1 \, $, and
$(u - \tilde{u}|\phi)_{\nabla} = 0$ for all $\phi \in \mathring{W}_2^1 \, $, forcing $ \, u - \tilde{u} = 0$.  \textbf{Done.}
\section{Extensions}
We introduce notation.  Write $\{ u=S(f,g) \}$ to express dependence of solution ($u \in W_2^1$) for (\ref{Diri-prob}) and (\ref{Diri-eqn}), on data $(f,g) \in \mathcal{L}_2 \times W_2^1 \,$.
What was denoted before by $\Lambda \phi$, denote now to show dependence on data: $\Lambda (\phi \,;\,f,g)$.
Let $J$ denote the duality-map on $ \mathring{W}_2^1 \,$, which maps $\phi$ to the function $(\phi | \cdot)_{\nabla}$ in the dual-space of $ \mathring{W}_2^1 \,$;
write $J^{-1} \, = \, R$, Riesz-representation.

Our proof of Dirichlet's principle had $ \{ \, R \Lambda ( \cdot \,;\, f,g) \, = \, p \, \}$ and $ \{ \, u \, = \, p \, + \, g \}$, so the solution-map is
\begin{equation}
S(f,g) \; = \; R \Lambda ( \cdot \,;\, f,g) \; + \; g \; . \label{Diri-prin2}
\end{equation}

To the standard formulation of Dirichlet's principle (stated above), the following fact deserves to be added. \\
\textbf{Note.}  Solution-map $S$ is linear and continuous. \\
\textbf{Proof.}  Straightforward calculation shows $S$ is linear.  For continuity, we will compute a bound on $S$.
Definition of $\| \cdot \|_{1,2} \,$, and (\ref{Poin/F}), give
\begin{equation*}
\| \, R \Lambda ( \cdot \,;\, f,g)  \, \|_{1,2}^2 \; =
\; \| \, R \Lambda ( \cdot \,;\, f,g) \, \|_2^2 \; + \; \| \, R \Lambda ( \cdot \,;\, f,g) \, \|_{\nabla}^2 \; \leq \;
(a^2 \; + \; 1) \, \| \, R \Lambda ( \cdot \,;\, f,g) \, \|_{\nabla}^2 \; .
\end{equation*}
Recall $R$ preserves norm, and use (\ref{Diri-prin}):
\begin{equation}
\| \, R \big( \Lambda ( \cdot \,;\, f,g) \big) \, \|_{1,2}^2 \; \leq \; (a^2 +1) \, \| \Lambda ( \cdot \,;\, f,g) \, \|^2 \; \leq \;
(a^2 +1) \, \Big( a \|f\|_2 \; + \; \|g\|_{\nabla} \Big)^2 \; .  \label{Diri-prin3}
\end{equation}
Use (\ref{Diri-prin3}) in (\ref{Diri-prin2}):
\begin{equation}
\| S(f,g) \|_{1,2} \; \leq \; \sqrt{a^2 + 1} \, \Big( a \|f\|_2 \; + \; \|g\|_{\nabla} \Big) \; + \; \| g \|_{1,2} \; . \label{Diri-prin4}
\end{equation}
(\ref{Diri-prin4}) implies continuity of $S$.  \textbf{Done.}

Another (expected) fact deserving addition to Dirichlet's principle is that $S(f,g)$ does not vary with choice $g$ for initial-extension of border-data $g_0 \,$.
As usual, we will express this without explicit use of $g_0 \,$. \\
\textbf{Note.}  $S(f,g) \; = \; S(f,\tilde{g}) \,$, ~if $\;f \in \mathcal{L}_2 \,$, $\; \{ g, \tilde{g} \} \subset W^1_2 \,$, $\; g - \tilde{g} \in \mathring{W}_2^1 \;$. \\
\textbf{Proof.}  Because $S$ is linear, it is enough to prove $\{ S(0,\psi)=0 \}$ if $\, \psi \in \mathring{W}^1_2 \;$.  Do this using (\ref{Diri-prin2}):
\begin{equation*}
S(0,\psi) \; = \; R \Lambda(\cdot \,;\, 0,\psi) \; + \; \psi \; = \; R \big( - (\psi| \cdot)_{\nabla} \big) \; + \; \psi \; = \; -\psi \; + \; \psi \; = \; 0 \; .
\end{equation*}
\textbf{Done.}

The preceding \textit{Note} says $S(f,g)$ depends on $g$ merely through its equi-class in quotient-space $ \, W_2^1 \Big/ \mathring{W}_2^1 \,$.
Hence the solution-map factors by a quotient-map; to help state this clearly, recall the following general fact. \\
\textbf{Note (A).}  Suppose $\{ A:X \to Y \}$ is a linear map between norm-spaces, $N$ is a closed subspace of $X$, $\{ q:X \to X/N \}$ is the quotient-map, and $Ax=0$ if $x \in N$.
Then there is a unique map $\{ \hat{A}: X/N \to Y \}$ with $A=\hat{A} \circ q \, $.  $ \hat{A} $ is linear.  $\hat{A}$ is continuous if $A$ is so. $\quad \Box$ \\

A factor of the solution-map will be the following map.
$$ Q \, : ~ \mathcal{L}_2 \times  W^1_2 \; \owns \, (f,g) \; \mapsto \; (f,g+\mathring{W}^1_2) \; \in \; \mathcal{L}_2 \times W^1_2 \Big/ \mathring{W}^1_2 \; .$$
$Q$ identifies naturally with a quotient-map, as noted in the following proof.  Now we are prepared to factor the solution-map. \\
\textbf{Note.}  There is a unique map $\{ \hat{S}: \mathcal{L}_2 \times W^1_2 \Big/ \mathring{W}^1_2 \to W^1_2 \} \,$ with $ \, S= \hat{S} \circ Q$.  $\hat{S}$ is linear and continuous. \\
\textbf{Proof.}  It is easy to check the following map is an isomorphism.
$$ \Phi \, : ~ \mathcal{L}_2 \times W^1_2 \Big/ \mathring{W}^1_2 \; \owns \, (f,g+\mathring{W}^1_2) \; \mapsto
\; (f,g) \; + \; \big(0 \times \mathring{W}^1_2 \big) \; \in \; \big( \mathcal{L}_2 \times  W^1_2 \big) \Big/ \big( 0 \times \mathring{W}^1_2 \big) \; .$$
Apply \textit{Note (A)} with $X=\mathcal{L}_2 \times W^1_2 \,$, $Y=W^1_2 \,$, $N=0 \times \mathring{W}^1_2 \,$, $A=S$; get $S=\hat{A} \circ q \,$.
Observe $\Phi \circ Q = q \,$; put $\hat{A} \circ \Phi = \hat{S}$; then $S=\hat{A} \circ q = \hat{A} \circ \Phi \circ Q = \hat{S} \circ Q$.  \textbf{Done. \\}
Call $\hat{S}$ the ``quotient solution-map''.  It will help express the solution neatly, in the following special case where data for the problem appears as a border-function.

Instead of our ``abstract'' quotient-space $\, W^1_2 (\Omega) \Big/ \mathring{W}^1_2 (\Omega) \,$, classical treatments of Dirichlet's problem used ``concrete'' spaces of functions on the border of region $\Omega$.
With our choice of space $\, W^1_2 (\Omega) \,$, concrete treatment requires more premises about region $\Omega$, to ensure existence of a (unique, bound-linear) trace-map $T$.
(Suitable premises are given in [A], [Z2], \textit{et al}.)
In that case, recall $T$ maps $W^1_2 (\Omega)$ \textit{onto} a Sobolev space for the border: $W^{1\!/2}_2 (\text{bdr} \, \Omega)\,$.
Null-space of $T$ is $\mathring{W}^1_2 (\Omega)\,$;
hence the quotient trace-map $\hat{T}$ is an isomorphism (of norm-spaces), from $W^1_2(\Omega) \big/ \mathring{W}^1_2(\Omega)$ to $W^{1\!/2}_2 (\text{bdr} \, \Omega)\,$;
so, in the present favorable case, we succeed to identify our abstract space with a concrete space of border-functions.

Let $E$ denote the inverse isomorphism to $\hat{T}$; border-data $g_0$ maps to its ``extension-class'' $E(g_0)$, consisting of those functions on $\Omega$ (in $W^1_2$) that extend $g_0 \,$ on $\mathrm{bdr}(\Omega)$.
If we follow the common practice of fixing datum $f \in \mathcal{L}_2 \,$, then Dirichlet's problem appears (as usual) as a border-data problem, whose solution-map is
\begin{equation}
W^{1\!/2}_2 (\text{bdr} \, \Omega)\; \owns \; g_0 \; \mapsto \; \hat{S}\big( \, f \, , \, E(g_0) \, \big) \; \in \; W^1_2(\Omega) \; . \label{Diri-prin5}
\end{equation}
In (\ref{Diri-prin5}), dependence of solution on border-data $g_0$ is affine (linear $+$ constant) and continuous.
Dependence is linear, in the important case ($f=0$).
Thus we neatly recover standard conclusions of classical treatments:
existence, uniqueness, and stability of solution for Dirichlet's border-data problem.

\textbf{References} \\
$\mathrm{[A]} \ \ \ $Adams, R.A.$\ \ \ $\textit{Sobolev Spaces}.  Academic Press, 1975. \\
$\mathrm{[G]} \ \ \ $G\aa rding, L.$\ \ \ $The Dirichlet Problem. \textit{The Mathematical Intelligencer}, 2(1), p43-53.  1979. \\
$\mathrm{[Z]} \ \ \ $Zeidler, E.$\ \ $\textit{Applied Functional Analysis: Applications to Mathematical Physics}.  Springer Pub., 1995. \\
$\mathrm{[Z2]} \ \ \ $Zeidler, E.$\ \ $\textit{Nonlinear Functional Analysis and its Applications}, vol. II/B.  Springer Pub., 1990.
\end{document}